\documentclass[11pt]{amsart}

\DeclareMathSymbol{\twoheadrightarrow} {\mathrel}{AMSa}{"10}

\def\Q{{\mathbf Q}}
\def\Z{{\mathbf Z}}
\def\C{{\mathbf C}}

\def\F{{\mathbf F}}
\def\L{{\mathbf L}}

\def\A{{\mathbf A}}
\def\Sn{{\mathbf S}_n}
\def\An{{\mathbf A}_n}

\def\Gal{\mathrm{Gal}}

\def\AGL{\mathrm{AGL}}
\def\Perm{\mathrm{Perm}}

\def\ord{\mathrm{ord}}
\def\Sz{\mathrm{Sz}}

\def\UU{\mathrm{U}}
\def\Fr{\mathrm{Fr}}
\def\pr{\mathrm{pr}}

\def\End{\mathrm{End}}

\def\Aut{\mathrm{Aut}}
\def\Hom{\mathrm{Hom}}

    \def\RR{\mathfrak{R}}

\def\fchar{\mathrm{char}}

\def\GL{\mathrm{GL}}
\def\sign{\mathrm{sign}}
\def\PGL{\mathrm{PGL}}
\def\PSL{\mathrm{PSL}}
\def\PSU{\mathrm{PSU}}

\def\SL{\mathrm{SL}}

\def\M{\mathrm{M}}

\def\dim{\mathrm{dim}}

\def\P{{\mathbf P}}

\newtheorem{thm}{Theorem}[section]
\newtheorem{lem}[thm]{Lemma}
\newtheorem{cor}[thm]{Corollary}
\newtheorem{prop}[thm]{Proposition}
\theoremstyle{definition}
\newtheorem{defn}[thm]{Definition}
\newtheorem{ex}[thm]{Example}
\newtheorem{exs}[thm]{Examples}

\newtheorem{rem}[thm]{Remark}

\title[Homomorphisms of hyperelliptic jacobians]
{Homomorphisms of hyperelliptic jacobians}
\author[Yu. G. Zarhin]{Yu. G. Zarhin}
\thanks{Partially supported by the NSF}

\begin{document}

\maketitle
\section{Definitions, notations, statements}

Let $K$ be a field. Let us fix its algebraic closure $K_a$ and
denote by $\Gal(K)$ the absolute Galois group  $\Aut(K_a/K)$ of
$K$. If $X$ is an abelian variety over $K_a$ then we write
$\End(X)$ for the ring of all its $K_a$-endomorphisms. If $Y$ is
(may be another) abelian variety over  $K_a$ then we write
$\Hom(X,Y)$ for the group of all $K_a$-homomorphisms from  $X$ to
$Y$. It is well-known that $\Hom(X,Y)=0$ if and only if
$\Hom(Y,X)=0$. One may easily check that if  $\End(X)=\Z$ and
$\dim(X)\ge\dim(Y)$ then $\Hom(X,Y)=0$ if and only if $X$ and $Y$
are {\sl not} isogenous over $K_a$.

Let $f(x) \in K[x]$ be a polynomial of degree $n\ge 3$ without
multiple roots.  We write $\RR_f\subset K_a$ for the set of its
roots,  $K(\RR_f) \subset K_a$ for the splitting field of $f$ and
$\Gal(f)=\Aut(K(\RR_f)/K)=\Gal(K(\RR_f)/K)$ for the Galois group
of $f$. It is well-known that $\RR_f$ consists of $n=\deg(f)$
elements. The group $\Gal(f)$ permutes elements of $\RR_f$ and
therefore can be identified with a certain subgroup of the group
$\Perm(\RR_f)$ of all permutations of  $\RR_f$. Clearly, every
ordering of  $\RR_f$ provides an isomorphism between
$\Perm(\RR_f)$ and the full symmetric group $\Sn$ which makes
$\Gal(f)$ a certain subgroup of $\Sn$. (It is well-known that this
permutation subgroup is transitive if and only if $f$ is
irreducible over $K$.)

Let us assume that $\fchar(K)\ne 2$ and consider the hyperelliptic
curve
$$C_f: y^2=f(x),$$
defined over $K$. Its genus $g=g(C_f)$ equals $(n-1)/2$ if $n$ is
odd and $(n-2)/2$ if  $n$ is even. Let $J(C_f)$ be the jacobian of
$C_f$; it is a $g$-dimensional abelian variety over $K_a$ that is
defined over $K$.

In his previous papers ~\cite{ZarhinMRL,ZarhinMRL2,ZarhinMMJ} the
author proved the following assertion.

\begin{thm}
\label{endo} Let $K$ be a field of characteristic different from
$2$. Let $n\ge 5$ be a positive integer. Let$f(x)\in K[x]$ be an
irreducible polynomial of degree $n \ge 5$. Assume also that if
 $\fchar(K)>0$ then $n \ge 9$ and
 $f(x)$ has no multiple roots. Suppose that the Galois group
of $f(x)$ coincides either with the full symmetric group $\Sn$ or
with the alternating group $\An$.

 Then  $\End(J(C_f))=\Z$.
\end{thm}

The main result of the present paper is the following statement.

\begin{thm}[Main Theorem]
\label{main} Let $K$ be a field of characteristic different from
$2$ and $K_a$ its algebraic closure. Let $f(x), h(x)\in K[x]$ be
irreducible polynomials of degree $n\ge 3$ and $m\ge 3$
respectively. Suppose that the splitting fields of
 $f$ and $h$ are linearly disjoint over  $K$.
Assume also that if $\fchar(K)>0$ then $n=\deg(f) \ge 9$ and
$f(x)$ and  $h(x)$ have no multiple roots.

Suppose that the following conditions hold:
\begin{enumerate}
\item[(i)]
$\Gal(h)=\A_m$ or $\mathbf{S}_m$.
\item[(ii)]
Either $\Gal(f)=\Sn$ or $\Gal(f)=\An$ and $n \ge 5$.
\end{enumerate}
Then
$$\Hom(J(C_f), J(C_h))=0,\quad \Hom(J(C_h),J(C_f))=0.$$
\end{thm}

We prove Theorem \ref{main} in \S \ref{proofm}.

\begin{ex}
Let $n\ge 3$ be a positive integer. It is well-known
  ~\cite[p. 139]{SerreM} that the Galois group of the polynomial $x^n-x-t$
  over the field of rational functions
$\Q(t)$ coincides with the full symmetric group $\Sn$. It follows
from Hilbert's irreducibility theorem that there exists an
 {\sl infinite}
set of rational numbers  $S\subset \Q$ such that for each $r\in S$
the Galois group  $\Gal(u_r)$ of the polynomial
$$u_r(x)=x^n-x-r\in \Q[x]$$ coincides with $\Sn$ and for {\sl
distinct} $r,k\in S$  the splitting fields of $u_r$ and $u_k$ are
linearly disjoint over $\Q$. Let us consider the jacobians
$J(C_{u_r})$ and $J(C_{u_k})$ of the hyperelliptic curves
$C_{u_r}:y^2=u_r(x)$ and $C_{u_k}:y^2=u_k(x)$ defined over $\Q$.
Notice that if $n<5$ then $J(C_{u_r})$ and $J(C_{u_k})$ are
elliptic curves. Applying Theorems  \ref{main} and \ref{endo} to
$u_r$ and $u_k$, we obtain that the jacobians $J(C_{u_r})$ are
$J(C_{u_k})$ absolutely simple and mutually non-isogenous over
$\bar{\Q}$ (and therefore over $\C$) for all
 $n \ge 3$ . In particular, for each positive integer
 $g$ the set of isogeny classes of absolutely simple
$g$-dimensional abelian varieties over $\Q$ is infinite. (This
assertion is well-known in the case of elliptic curves.) It also
follows from Theorem \ref{endo} that for each positive integer
 $g>1$ the set of isogeny classes of absolutely simple
$g$-dimensional abelian varieties over $\Q$ without nontrivial
endomorphisms over $\C$ is infinite. (The similar assertion for
elliptic curves is also well-known: it suffices to take for each
prime $p$ an elliptic curve with $j$-invariant $1/p$.)
\end{ex}

\begin{cor}
\label{NneM} Let $K$ be a field of characteristic different from
$2$ and $K_a$ its algebraic closure. Let $f(x), h(x)\in K[x]$ be
irreducible polynomials of degree $n\ge 5$ and $m\ge 3$
respectively. Assume also that if $\fchar(K)>0$ then $n=\deg(f)
\ge 9$ and both polynomials $f(x)$ and $h(x)$ have no multiple
roots. Suppose that the following conditions hold:

\begin{enumerate}
\item[(i)]
$\Gal(h)=\A_m$ or $\mathbf{S}_m$.
\item[(ii)] $\Gal(f)=\Sn$ or $\Gal(f)=\An$.
\item[(iii)]
Either $n \ne m$ or  $\Gal(f)=\Sn$ and $\Gal(h)=\A_m$.
\end{enumerate}
Then
$$\Hom(J(C_f), J(C_h))=0, \quad \Hom(J(C_h),J(C_f))=0.$$
\end{cor}

\begin{cor}
\label{NeqM} Let $K$ be a field of characteristic different from
$2$ and $K_a$ its algebraic closure.  Let $n \ge 5$ be a positive
integer different from $6$. Let $f(x), h(x)\in K[x]$ be
irreducible polynomials of degree $n$. Assume also that if
$\fchar(K)>0$ then $n \ge 9$ and both polynomials $f(x)$ and
$h(x)$ have no multiple roots. Suppose that the following
conditions hold:

\begin{enumerate}
\item[(i)]
$\Gal(h)=\An$ or $\Sn$.
\item[(ii)] $\Gal(f)=\Sn$ or $\Gal(f)=\An$.
\item[(iii)]
Let us put $K_f:=K[x]/f K[x], K_h:=K[x]/h K[x]$. Then the field
extensions $K_f/K$ and $K_h/K$ are not isomorphic.
\end{enumerate}
Then
$$\Hom(J(C_f), J(C_h))=0, \quad \Hom(J(C_h),J(C_f))=0.$$
\end{cor}

We will prove Corollaries \ref{NneM} and \ref{NeqM} in \S
\ref{proofc}.

\section{Proof of Main Theorem}
\label{proofm} Let $d$ be a positive integer that is not divisible
by $\fchar(K)$.  Let $X$ be an abelian variety of positive
dimension defined over $K$. We write $X_d$ for the kernel of
multiplication by $d$ in $X(K_a)$. It is known \cite{Mumford} that
the commutative group  $X_d$ is a free $\Z/d\Z$-module of rank
$2\dim(X)$. Clearly,  $X_d$ is a Galois submodule in  $X(K_a)$. We
write
$$\tilde{\rho}_{d,X}:\Gal(K) \to \Aut_{\Z/d\Z}(X_d) \cong
\GL(2\dim(X),\Z/d\Z)$$ for the corresponding (continuous)
homomorphism defining the Galois action on $X_d$. Let us put
$$\tilde{G}_{d,X}=\tilde{\rho}_{d,X}(\Gal(K)) \subset
\Aut_{\Z/d\Z}(X_d).$$ Clearly,  $\tilde{G}_{d,X}$ coincides with
the Galois group of the field extension  $K(X_d)/K$ where $K(X_d)$
is the field of definition of all points of order dividing $d$ on
$X$. In particular, if
 $\ell\ne \fchar(K)$ is a prime then
$X_{\ell}$ is a $2\dim(X)$-dimensional vector space over the prime
field $\F_{\ell}=\Z/\ell\Z$ and the inclusion $\tilde{G}_{\ell,X}
\subset \Aut_{\F_{\ell}}(X_{\ell})$ defines a faithful linear
representation of the group  $\tilde{G}_{\ell,X}$ in the vector
space $X_{\ell}$.
 We will deduce Theorem \ref{main} from the following auxiliary statement which is of some independent
 interest.

\begin{thm}
 \label{mainA}
 Let $\ell$ be a prime, $K$ a field of characteristic different from
 $\ell$, $X$ and $Y$  are abelian varieties of positive dimension defined over $K$.
 Suppose that the following conditions hold:

\begin{enumerate}
\item[(i)]
The extensions $K(X_{\ell})$ and $K(Y_{\ell})$  are linearly
disjoint over $K$.
\item[(ii)]
The natural representation of the group $\tilde{G}_{\ell,X}$ in
$X_{\ell}$ is absolutely irreducible.
\item[(iii)]
The natural representation of the group $\tilde{G}_{\ell,Y}$ in
$Y_{\ell}$ is irreducible.
\end{enumerate}

Then either
$$\Hom(X,Y)=0, \quad \Hom(Y,X)=0$$
or $\fchar(K)>0$ and both abelian varieties $X$ and $Y$ are
supersingular.
\end{thm}

We will prove Theorem \ref{mainA}  in   \S \ref{proofA}.

In fact, we are going to prove not Theorem \ref{mainA} but its
certain generalization. In order to state this generalization, we
need to introduce definitions of {\sl nice} and {\sl very nice}
polynomials. But first, let us recall some standard notations
~\cite[\S 2.8]{DM}. Hereafter, $\F_q$ denotes the $q$-element
finite field of characteristic $p$, $\GL(d,q):=\GL(d,\F_q)$
denotes the group of invertible linear transformations of the
$d$-dimensional vector space $\F_q^d$, $\SL(d,q):=\SL(d,\F_q)$ its
subgroup of all matrices with determinant $1$ and
$\PGL(d,q)=\PGL(d,\F_q)$ and $\L_d(q)=\PSL(d,q)=\PSL(d,\F_q)$ are
the corresponding quotients with respect the subgroups of scalar
matrices, viewed as transformation groups of the projective space
$\P^{d-1}(\F_q)$. In addition,
 $\AGL(d,q):=\AGL(d,\F_q)$ is the group of all affine transformations of  $\F_q^d$,
which is a semi-direct product of
 $\GL(d,q)$ and the group $\F_q^d$ of all translations.
We write $\Fr:\F_q^d\to \F_q^d$ for the Frobenius automorphism
$$(a_1,\cdots , a_d) \mapsto (a_1^p,\cdots , a_d^p).$$ We write
$\Gamma\L(d,q)$, $\Sigma\L(d,q)$ and $\A\Sigma\L(d,q)$ for the
transformation groups of $\F_q^d$, generated by $\Fr$ and
 $\GL(d,q),\SL(d,q)$ and $\AGL(d,q)$ respectively.
We write $\P\Gamma\L(d,q)$ and $\P\Sigma\L(d,q)$ for the
transformation groups of $\P^{d-1}(\F_q)$ induced by
$\Gamma\L(d,q)$ and $\Sigma\L(d,q)$ respectively. (In other words,
$\P\Gamma\L(d,q)$ and $\P\Sigma\L(d,q)$ are the quotients of
$\Gamma\L(d,q)$ and $\Sigma\L(d,q)$ respectively with respect to
the corresponding subgroups of scalar matrices.)

 Let $f(x) \in K[x]$ be a separable irreducible
polynomial of degree  $n \ge 3$. We say that $f$ is {\sl very
nice} if one of the following conditions holds:
\begin{enumerate}
\item[(s)]
$\Gal(f)=\Sn$.
\item[(a)]
$\Gal(f)=\An$ and $n \ge 5$.
\item[(m)]
 $n=11$ or $12$ and $\Gal(f)$ is the corresponding small Mathieu group
$\M_n$  acting $4$(or $5$-)transitivelyon $\RR_f$.
\item[(l11)]
$n=11$ and $\Gal(f)=\L_2(11)=\PSL_2(\F_{11})$ acts doubly
transitively on  $\RR_f$.
\item[(m12)]
$n=12$ and $\Gal(f)=\M_{11}$ acts $3$-transitively  on $\RR_f$.
\item[(aff)]
There exist an odd prime   $p$, its positive integral power  $q$
and a positive integer $d$ such that $n=p^d>3$ and one may
identify $\RR_f$ with $\F_q^d$ in such a way that $\Gal(f)$
becomes {\sl $2$ or $3$-transitive} subgroup of $\AGL(d,q)$ that
contains the subgroup $\F_q^d$ of all translations.
\item[(p)]
There exist an odd prime   $p$, its positive integral power $q$
and a positive integer  $d\ge 3$ such that $n=\frac{q^d-1}{q-1}$
and one may identify $\RR_f$ with  $\P^{d-1}(\F_q)$ in such a way
that $\Gal(f)$ becomes a subgroup of $\P\Gamma\L(d,q)$ that
contains $\PSL(d,q)$.
\item[(p1)]
There exist an odd prime  $p$ and its positive integral power $q$
such that $n=q+1$ and one may identify
  $\RR_f$ with the projective line $\P^{1}(\F_q)$ in such a way
that $\Gal(f)$ becomes a $3$-transitive subgroup of
$\P\Gamma\L(2,q)$.
\item[(p2)]
There exists a positive integer $d\ge 2$ such that  $q:=2^d,n=q+1$
and one may identify $\RR_f$ with the projective line
$\P^{1}(\F_q)$ in such a way that  $\Gal(f)$ becomes a subgroup of
$\P\Gamma\L(d,q)$ that contains $\PSL(2,q)$.
\item[(u3)]
There exists a positive integer $d\ge 2$ such that
$q:=2^d,n=q^3+1$, and one may identify  $\RR_f$ with the set of
isotropic lines  (Hermitian curve) in $\F_{q^2}^3$ with respect to
a certain non-degenerate Hermitian form in such a way that
 $\Gal(f)$ becomes a group that contains the corresponding projective special
 unitary group
$\UU_3(q):=\PSU(3,q)=\PSU(3,\F_{q^2})$ and  $\UU_3(q)$ acts doubly
transitively on $\RR_f$.
\item[(sz)]
There exists a positive integer  $d$ such that $q:=2^{2d+1},
n=q^2+1$ and  $\Gal(f)$ contains a subgroup isomorphic to the
Suzuki group $\Sz(q)$ and $\Sz(q)$ acts doubly transitively on
$\RR_f$.
\end{enumerate}
A polynomial $f$ is called {\sl nice} if either it is very nice or
one of the following conditions holds:
\begin{enumerate}
\item[(a3)]
$n=3$ and $\Gal(f)=\A_3$.
\item[(a4)]
$n=4$ and $\Gal(f)=\A_4$.
\item[(p3)]
There exist an odd prime  $p$ and its positive integral power $q$
such that $n=q+1$, and one may identify $\RR_f$  with the
projective line $\P^{1}(\F_q)$ in such a way that $\Gal(f)$
becomes a doubly transitive subgroup of $\P\Gamma\L(2,q)$. In
addition, $q$ must be congruent either to  $3$ or to $5$ modulo
$8$.
\end{enumerate}

\begin{rem}
The doubly transitive action of the Suzuki groups
 $\Sz(q)$
(the case (sz)) is described explicitly on pp.  184--187 of
\cite{HB}; see \cite{ZarhinTexel} concerning  the relations to
hyperelliptic jacobians. Concerning the doubly transitive action
of $\UU_3$ on the Hermitian curve (the case (u3)) see ~\cite[Kap.
II, Satz 4.12]{H}, ~\cite[pp. 248--250]{DM}; the relations with
hyperelliptic jacobians are discussed in \cite{ZarhinPAMS}.
\end{rem}

In order to explain what nice polynomials are good for, let us
recall the definition of the   {\sl heart} of the permutational
action of $\Gal(f)$  on
 $\RR_f$ (\cite{Mortimer},
\cite{ZarhinTexel}).

 Let $\RR =\RR_f=\{a_1, \dots , a_n\}\subset K_a$ be the set of
all roots of $f$. We may view $\Sn$ as the group of all
permutations of  $\RR$. The Galois group $G=\Gal(f)$ of $f$
permutes the roots and therefore becomes a subgroup of $\Sn$. The
action of $G$ on $\RR$ defines the standard {{\sl permutational}
representation in the $n$-dimensional $\F_{2}$-vector space
$\F_2^{\RR}$ of all functions $\psi:R \to \F_{2}$. This
representation is not irreducible. Indeed, the "line" of constant
functions
 $\F_{2}\cdot 1$  and the hyperplane
$(\F_{2}^{\RR})^0:=\{\psi\mid \sum_{i=1}^n \psi(a_i)=0\}$ are
$G$-invariant subspaces in $\F_2^{\RR}$.
 If $n$ is odd then one calls  $(\F_{2}^{\RR})^0$  the {\sl heart}
 of the permutational action of
 $G=\Gal(f)$ on
$\RR=\RR_f$ over $\F_2$ and denotes it by $Q_{\RR}=Q_{\RR_f}$. If
$n$ is even then $(\F_{2}^{\RR})^0$ contains $\F_{2}\cdot 1$ and
we obtain the natural representation of  $G=\Gal(f)$ in the
$(n-2)$-dimensional $\F_{2}$-vector quotient-space
$$(\F_{2}^B)^{00}:=(\F_{2}^{\RR})^0/(\F_{2}\cdot 1).$$
In this case $(\F_{2}^B)^{00}$ is also called the {\sl heart}
 of the permutational action of
 $G=\Gal(f)$ on
$\RR=\RR_f$ over $\F_2$ and denoted by $Q_{\RR}=Q_{\RR_f}$.

It is known \cite{Klemm} that if $n$ is odd and the
$\Gal(f)$-module $Q_{\RR_f}$ is absolutely simple then $\Gal(f)$
acts on $\RR_f$ doubly transitively.

\begin{rem}
\label{goodirr}
 If a polynomial $f(x)$ is nice then:
 \begin{enumerate}
 \item[(i)]
 Either $(n,\Gal(f))=(3,\A_3)$ or $\Gal(f)$ acts doubly transitively
 on
 $\RR_f$;
\item[(ii)]
The $\Gal(f)$-module $Q_{\RR_f}$ is simple. In addition,
$Q_{\RR_f}$ is absolutely simple if and only if  $f(x)$ is {\sl
very} nice. In the case of doubly transitive $\Gal(f)$ this
assertion follows immediately from results of
~\cite{Mortimer,Ivanov}. The remaining case  $n=3, \Gal(f)=\A_3$
 is easy. (See also ~\cite{ZarhinTexel,ZarhinPAMS}.)
\end{enumerate}
\end{rem}
\begin{rem}
Let us assume that a permutation group  $\Gal(f)\subset
\Perm(\RR_f)$ is isomorphic to one of the {\sl known doubly
transitive} permutation groups ~\cite[\S 7.7]{DM}. Then $f(x)$ is
nice if and only if the $\Gal(f)$-module $Q_{\RR_f}$ is simple.
This assertion follows easily from results of
~\cite{Mortimer,Ivanov}.
\end{rem}

Now we are ready to state the promised generalization of Main
Theorem.

\begin{thm}
\label{mgood} Let $K$ be a field of characteristic different from
$2$ and $K_a$ its algebraic closure.  Let $f(x), h(x)\in K[x]$ be
irreducible polynomials without multiple roots of degree
 $n\ge 3$ and
$m\ge 3$ respectively. Suppose that the splitting fields of  $f$
and $h$ are linearly disjoint over $K$. Suppose that $f(x)$ is
very nice and $h(x)$ is nice.

Then either
$$\Hom(J(C_f), J(C_h))=0, \quad \Hom(J(C_h),J(C_f))=0$$
or $\fchar(K)>0$ and both jacobians $J(C_f)$ and  $J(C_h)$ are
supersingular abelian varieties.
\end{thm}

\begin{proof}[Proof of Theorem \ref{mgood}]
The canonical surjection $\Gal(K) \twoheadrightarrow \Gal(f)$
defines on the $\Gal(f)$-module $Q_{\RR_f}$ the natural structure
of $\Gal(K)$-module. It is well-known that the  $\Gal(K)$-modules
$Q_{\RR_f}$ and $J(C_f)_2$ are canonically isomorphic (see, for
instance, \cite{MumfordENS}, \cite{Mori2} or \cite{ZarhinTexel}).
This implies, in particular, in light of Remark
  \ref{goodirr},
that the $\tilde{G}_{2,J(C_f)}$-module is absolutely simple.
Similarly, the canonical surjection  $\Gal(K) \twoheadrightarrow
\Gal(h)$ provides the $\Gal(h)$-module $Q_{\RR_h}$ with natural
structure of the $\Gal(K)$-module and the  $\Gal(K)$-modules
$Q_{\RR_h}$ and $J(C_h)_2$ are canonically isomorphic. Now it
follows from Remark  \ref{goodirr} that the
$\tilde{G}_{2,J(C_h)}$-module is simple.  We have
$$K(J(C_f)_2)\subset K(\RR_f),\quad K(J(C_h)_2)\subset K(\RR_h).$$
Since the field extensions $K(\RR_f)/K$ and $K(\RR_h)/K$ are
linearly disjoint, their subextensions $K(J(C_f)_2)/K$ and
$K(J(C_h)_2)/K$ are also linearly disjoint. One has only to apply
Theorem  \ref{mainA} to $\ell=2, X=J(C_f), Y=J(C_h)$.
\end{proof}

\begin{rem}
In fact, if $n \ne 4$ (respectfully $m \ne 4$) then the
$\Gal(f)$-module $Q_{\RR_f}$ is faithful and $K(J(C_f)_2)=
K(\RR_f)$ (respectfully the $\Gal(h)$-module $Q_{\RR_h}$ is
faithful and $K(J(C_h)_2)= K(\RR_h)$).
\end{rem}

\begin{proof}[Proof of Theorem \ref{main}]
It follows from Theorem  \ref{mgood} that if there exists a
non-zero homomorphism between  $J(C_f)$ and $J(C_h)$ then
$\fchar(K)>0$ and both jacobians are supersingular. However if
 $\fchar(K)>0$ then $n \ge 9$ and,  thanks to Theorem
\ref{endo}, $\End(J(C_f)=\Z$ and therefore $J(C_f)$ is {\sl not}
supersingular.
\end{proof}

\section{Homomorphisms of abelian varieties}
\label{proofA} In order to prove Theorem \ref{mainA}, we need the
following elementary statement that is well known when the ground
field is algebraically closed and has characteristic zero
(~\cite[\S 3.2]{Serre}; see also theorem 10.38 of \cite{CR}).

\begin{lem}
\label{irrprod} Let $F$ be a field. Let $H_1$ and $H_2$ be groups.
Let  $\tau_1: H_1 \to \Aut_F(W_1)$ be an irreducible
finite-dimensional representation of  $H_1$ over $F$ and $\tau_2:
H_2 \to \Aut_F(W_2)$  be an absolutely irreducible
finite-dimensional representation of  $H_2$ over $F$. Then the
natural linear representation

$$\tau_1^*\otimes\tau_2:H_1\times H_2 \to \Aut_F(\Hom_F(W_1,W_2))$$

of the group $H_1\times H_2$ in the $F$-vector space
$\Hom_F(W_1,W_2)$ is irreducible.
\end{lem}

\begin{rem}
Clearly, the representations of $H_1\times H_2$ in
$\Hom_F(W_1,W_2)$ and $\Hom_F(W_2,W_1)$ are mutually dual.
Therefore the irreducibility of  $\Hom_F(W_1,W_2)$ implies the
irreducibility of  $\Hom_F(W_2,W_1)$.
\end{rem}

We will prove Lemma \ref{irrprod} at the end of this Section.

\begin{proof}[Proof of Theorem \ref{mainA}]
First, notice that the natural representation
$$\Gal(K)
\to\Aut_{\F_{\ell}}(\Hom_{\F_{\ell}}(Y_{\ell},X_{\ell}))$$ is
irreducible. Indeed, let us denote this representation by $\tau$
and let us put
$$F=\F_{\ell}, H_1=\tilde{G}_{\ell,Y}, W_1=Y_{\ell},
H_2=\tilde{G}_{\ell,X}, W_2=X_{\ell}.$$ Denote by
$$\tau_1:H_1=\tilde{G}_{\ell,Y}\subset
\Aut_{\F_{\ell}}(Y_{\ell})=\Aut_{\F_{\ell}}(W_1)$$ and
$$\tau_2:H_2=\tilde{G}_{\ell,X}\subset
\Aut_{\F_{\ell}}(X_{\ell})=\Aut_{\F_{\ell}}(W_2)$$ the
corresponding inclusion maps.

It follows from Lemma  \ref{irrprod} that the linear
representation
$$\tau_1^*\otimes\tau_2:\Gal(K(Y_{\ell})/K)\times \Gal(K(X_{\ell})/K)
\to\Aut_{\F_{\ell}}(\Hom_{\F_{\ell}}(Y_{\ell},X_{\ell}))$$ is
irreducible.

One may easily check that the homomorphism  $\tau$, which defines
the structure of  $\Gal(K)$-module on
$\Hom_{\F_{\ell}}(Y_{\ell},X_{\ell})$, coincides with the
composition of the natural surjection  $\Gal(K) \twoheadrightarrow
\Gal(K(X_{\ell},Y_{\ell})/K)$, the natural embedding
$$\Gal(K(X_{\ell},Y_{\ell})/K)\hookrightarrow
\Gal(K(Y_{\ell})/K)\times \Gal(K(X_{\ell})/K)$$ and
$$\tau_1^*\otimes\tau_2:\Gal(K(Y_{\ell})/K)\times
\Gal(K(X_{\ell})/K)\to
\Aut_{\F_{\ell}}(\Hom_{\F_{\ell}}(Y_{\ell},X_{\ell})).$$ Here
 $K(X_{\ell},Y_{\ell})$ is the compositum of the fields
 $K(X_{\ell})$ and $K(Y_{\ell})$. The linear disjointness of
$K(X_{\ell})$ and $K(Y_{\ell})$ means that
$$\Gal(K(X_{\ell},Y_{\ell})/K)= \Gal(K(Y_{\ell})/K)\times
\Gal(K(X_{\ell})/K).$$ This implies that $\tau$ is the composition
of  {\sl surjective}  $\Gal(K) \twoheadrightarrow
\Gal(K(Y_{\ell})/K)\times \Gal(K(X_{\ell})/K)$ and
$\tau_1^*\otimes\tau_2$. Since the representation
$$\tau_1^*\otimes\tau_2:\Gal(K(Y_{\ell})/K)\times \Gal(K(X_{\ell})/K)
\to\Aut_{\F_{\ell}}(\Hom_{\F_{\ell}}(Y_{\ell},X_{\ell}))$$ is
irreducible, the representation $$\tau: \Gal(K)\to
\Aut_{\F_{\ell}}(\Hom_{\F_{\ell}}(Y_{\ell},X_{\ell}))$$ is also
irreducible.

Second, let $T_{\ell}(X)$ and $T_{\ell}(Y)$ be the Tate
$\Z_{\ell}$-modules of abelian varieties  $X$ and $Y$ respectively
\cite{Mumford}. Recall that $T_{\ell}(X)$ and $T_{\ell}(Y)$ are
free $\Z_{\ell}$-modules of rank $2\dim(X)$ and $2\dim(Y)$
respectively. There are also natural continuous homomorphisms
$$\rho_{\ell,X}:\Gal(K) \to \Aut_{\Z_{\ell}}(T_{\ell}(X)),\quad \rho_{\ell,Y}:\Gal(K) \to
\Aut_{\Z_{\ell}}(T_{\ell}(Y)).$$ There are also natural
homomorphisms
$$X_{\ell}=T_{\ell}(X)/\ell T_{\ell}(X), Y_{\ell}=T_{\ell}(Y)/\ell
T_{\ell}(Y),$$ which are isomorphisms of Galois modules. So one
may view $\tilde{\rho}_{\ell,X}$ as the reduction of
$\rho_{\ell,X}$ modulo $\ell$ and $\tilde{\rho}_{\ell,Y}$ as the
reduction of $\rho_{\ell,Y}$ modulo $\ell$. It is also convenient
to consider the Tate $\Q_{\ell}$-modules
$V_{\ell}(X)=T_{\ell}(X)\otimes_{\Z_{\ell}}\Q_{\ell}$ and
$V_{\ell}(Y)=T_{\ell}(Y)\otimes_{\Z_{\ell}}\Q_{\ell}$, which are
$\Q_{\ell}$-vector spaces of dimension   $2\dim(X)$ and $2\dim(Y)$
respectively. The groups  $T_{\ell}(X)$ and $T_{\ell}(Y)$ are
naturally identified with the $\Z_{\ell}$-lattices in
$V_{\ell}(X)$ and  $V_{\ell}(Y)$ respectively and the inclusions
$$\Aut_{\Z_{\ell}}(T_{\ell}(X))\subset
\Aut_{\Q_{\ell}}(V_{\ell}(X)), \quad
\Aut_{\Z_{\ell}}(T_{\ell}(Y))\subset
\Aut_{\Q_{\ell}}(V_{\ell}(Y))$$ allow us to consider $V_{\ell}(X)$
and $V_{\ell}(Y)$ as representations of $\Gal(K)$ over
$\Q_{\ell}$.

Third, I claim that the natural representation of  $\Gal(K)$ in
 $\Hom_{\Q_{\ell}}(V_{\ell}(X),V_{\ell}(Y))$
over $\Q_{\ell}$ is irreducible. Indeed, the $\Z_{\ell}$-module
$\Hom_{\Z_{\ell}}(T_{\ell}(X),T_{\ell}(Y))$ is a
$\Gal(K)$-invariant $\Z_{\ell}$-lattice in
$\Hom_{\Q_{\ell}}(V_{\ell}(Y),V_{\ell}(X))$. On the other hand,
the reduction of this lattice modulo $\ell$ coincides with
$$\Hom_{\Z_{\ell}}(T_{\ell}(Y),T_{\ell}(X))\otimes
\Z/\ell\Z=\Hom_{\F_{\ell}}(T_{\ell}(Y)/\ell
T_{\ell}(Y),T_{\ell}(X)/\ell
T_{\ell}(X))=\Hom_{\F_{\ell}}(Y_{\ell},X_{\ell}).$$ But we just
established the simplicity of the $\Gal(K)$-module
$\Hom_{\F_{\ell}}(Y_{\ell},X_{\ell})$. It follows easily that the
$\Gal(K)$-module $\Hom_{\Q_{\ell}}(V_{\ell}(X),V_{\ell}(Y))$ is
simple  (see, for instance, exercise  2 in \S 15.2 of Serre's book
\cite{Serre}).

Fourth, notice that there is a natural embedding  (\cite{Mumford},
\S 19)
$$\Hom(Y,X)\otimes\Q_{\ell}\subset\Hom_{\Q_{\ell}}(V_{\ell}(X),V_{\ell}(Y)),$$
whose image is a $\Gal(K)$-invariant subspace. The irreducibility
of $\Hom_{\Q_{\ell}}(V_{\ell}(X),V_{\ell}(Y))$ implies that either
$$\Hom(Y,X)\otimes\Q_{\ell}=\Hom_{\Q_{\ell}}(V_{\ell}(X),V_{\ell}(Y))$$
or $\Hom(Y,X)\otimes\Q_{\ell}=0$. Since $\Hom(Y,X)$ is a free
commutative group of finite rank, either $\Hom(Y,X)=0$ or the rank
of $\Hom(Y,X)$ equals $4\cdot\dim(X)\cdot\dim(Y)$. In order to
finish the proof, we need the following proposition.

\begin{prop}
\label{super} Let $A$ and $B$ are abelian varieties of positive
dimension over an algebraically closed field $\mathcal{K}$.
Suppose that the rank of the group $\Hom(A,B)$ equals
$4\cdot\dim(A)\cdot\dim(B)$. Then $\fchar(\mathcal{K})>0$ and both
$A$ and $B$ are supersingular.
\end{prop}

\begin{proof}[Proof of Proposition \ref{super}]
The case $A=B$ was treated in lemma 3.1 of  \cite{ZarhinMRL}.

 Replacing  $A$ and $B$ by isogenous abelian varieties, we may
 assume that they split into products

$$A=\prod_i A_i, \quad B=\prod_j B_j$$
of simple abelian varieties $A_i$ and $B_j$ respectively. Since

$$\dim(A)=\sum_i\dim(A_i),\dim(B)=\sum_j\dim(B_j),
\Hom(A,B)=\prod_{i,j}\Hom(A_i,B_j),$$ and the rank of the free
commutative group  $\Hom(A_i,B_j)$ does not exceed
$4\cdot\dim(A_i)\cdot\dim(B_j)$ (\cite{Mumford}, \S 19, corollary
1 to theorem 3), the rank of $\Hom(A_i,B_j)$ equals
$4\cdot\dim(A_i)\cdot\dim(B_j)$ for all $i$ and $j$. Since $A_i$
and $B_j$ are simple, they are isogenous. This implies that
$\dim(A_i)=\dim(B_j)$ and the rank of each of the free commutative
group (with respect to addition) $\End(A_i)$ and $\End(B_j)$
equals
$$4\cdot\dim(A_i)\cdot\dim(B_j)=4\cdot\dim(A_i)^2=4\cdot\dim(B_j)^2.$$
Applying lemma 3.1 of \cite{ZarhinMRL} to each $A_i$ and $B_j$, we
conclude that  $\fchar(\mathcal{K})>0$ and all $A_i$ and $B_j$ are
supersingular. It follows easily that $A$ and $B$ are also
supersingular.
\end{proof}
{\sl End of Proof of Theorem} \ref{mainA}. Applying Proposition
 \ref{super} to $A=Y$ and $B=X$, we conclude
that  $\fchar(K)>0$ and $X$ and $Y$ are supersingular.
\end{proof}

\begin{proof}[Proof of Lemma \ref{irrprod}]
Throughout the proof all the tensor products are taken over $F$.
 First, replacing the $H_1$-module $W_1$ by its dual
$W_1^*=\Hom_F(W_1,F)$, we reduce the problem to the assertion
about the irreducibility of the tensor product
$$\tau_1\otimes\tau_2:H_1\times H_2 \to \Aut_F(W_1\otimes
W_2).$$ Since the $H_2$-module $W_2$ is absolutely simple, the
corresponding $F$-algebra homomorphism
$$F[H_2] \to \End_F(W_2)$$
(induced by $\tau_2$) is surjective. Here $F[H_2]$ is the group
algebra of $H_2$.

Let us denote by  $D$ the endomorphism ring of the $F[H_1]$-module
$W_1$. Since $W_1$ is simple, $D$ is a division algebra, whose
center contains $F$. Clearly, the $F$-dimension of $D$ is finite
and $W_2$ is a free $D$-module of finite rank. It follows from
Jacobson's density theorem that the image of the (induced by
$\tau_2$)  $F$-algebra homomorphism
$$F[H_1] \to \End_F(W_1)$$
coincide with $\End_D(W_1)$. Here $F[H_1]$ is the group algebra of
 $H_1$. There is the natural structure of the free  $D\otimes F=D$-module of finite rank on
 $W_1\otimes W_2$. Clearly, the  $\End_D(W_1\otimes W_2)$-module $W_1\otimes W_2$
is simple.

It follows that the image of the (induced by
$\tau_1\otimes\tau_2$)  $F$-algebra homomorphism
$$F[H_1\times H_2]=F[H_1]\otimes F[H_2] \to \End_F(W_1\otimes W_2)$$
coincides with  $\End_D(W_1)\otimes \End_F(W_2)$. Applying lemma
10.37 on p.  252 of \cite{CR}, we conclude that
$$\End_D(W_1)\otimes \End_F(W_2)=\End_{D\otimes F}(W_1\otimes
W_2).$$ Therefore the image of the group algebra  $F[H_1\times
H_2]$ in $\End_F(W_1\otimes W_2)$ coincides with
$$\End_{D\otimes F}(W_1\otimes W_2)=\End_{D}(W_1\otimes W_2).$$
 Now the simplicity of the
$\End_D(W_1\otimes W_2)$-module $W_1\otimes W_2$ implies the
simplicity of the  $F[H_1\times H_2]$-module $W_1\otimes W_2$.
\end{proof}

\section{Proof of Corollaries \ref{NneM} and \ref{NeqM}}
\label{proofc} We start with the following useful definition.
\begin{defn}
Finite groups $G_1$ and $G_2$ are called {\sl disjoint} if they do
not have isomorphic quotients except the trivial one-element
group.
\end{defn}

\begin{exs}
\label{exd} Clearly, the following pairs provide examples of
disjoint groups.

\begin{enumerate}
\item[(i)]
$\mathbf{S}_3$ and $\A_3$;
\item[(ii)]
$\mathbf{S}_n$ and $\A_m$ ($m \ge 5$);
\item[(iii)]
$\An$ and $\A_m$ ($n\ne m$ and $m \ge 5$);
\item[(iv)]
$G_1:=\PSL(d,q)\subset G_2:=\PGL(d,q)$, where
\begin{enumerate}
\item[(a)]
 $d>1,(d,q) \ne (2,2), (d,q) \ne (2,3);$
 \item[(b)]
 integers $d$ and $q-1$ have a common divisor $>1$.
  \end{enumerate}
The  condition (a) means that $G_1$ is a finite simple non-abelian
group
 ~\cite[Ch. 1, \S 9]{Suzuki}. The  condition (b) means that $G_1\ne G_2$.
 Clearly, $G_1$ is a normal subgroup of   $G_2$ and the quotient  $G_2/G_1$ is a cyclic group of order
  $r$ where $r$ is the largest common divisor of  $d$ and
 $q-1$. In order to prove that $G_1$ and $G_2$ are disjoint, it
 suffices to check that there does not exist a surjective
 homomorphism
 $\phi: G_2 \twoheadrightarrow G_1$.
 Let us assume that such a surjection does exist.
Then its kernel
   $\ker(\phi)$ is a proper normal subgroup of
 $G_2=\PGL(d,q)$ and its preimage
 $G'$ in $\GL(d,q)$ is a proper normal subgroup of $\GL(d,q)$
 containing all the scalars and also an element that is not a
 scalar.
 Since every normal subgroup of
  $\GL(d,q)$ either contains $\SL(d,q)$ or consists of scalars
  \cite[Ch. 1, \S 9, Th. 9.9]{Suzuki}, we conclude that
  $G'$ contains $\SL(d,q)$ and therefore
 $\ker(\phi)$ contains $\PSL(d,q)=G_1$. This implies that the image
 $G_1$ of the surjection $\phi$ is isomorphic to a quotient of the
 cyclic group
  $G_2/G_1$ and therefore is also cyclic. Since $G_1:=\PSL(d,q)$
  is non-abelian, we obtain the desired contradiction, which
  proves the disjointness of
  $G_1$ and $G_2$.
\end{enumerate}

\end{exs}

Let us recall the statement of well-known Goursat's lemma
 (see for instance  ~\cite[p. 75]{Lang})

\begin{lem}
\label{goursat} Let $G_1$ and $G_2$ be finite groups. Let   $H$ be
a subgroup of the product $G_1\times G_2$ such that the
corresponding projection maps $\pr_1:H\to G_1$ and $\pr_2: H \to
G_2$ are surjective. Denote by $H_1$ (respectfully by $H_2$) the
normal subgroup of  $G_1$ (respectfully of $G_2$), such that the
kernel of  $\pr_2$ (respectfully of $\pr_1$) coincides with
$H_1\times \{1\}$ (respectfully with $\{1\}\times H_2$). Then
there exists an isomorphism $\gamma:G_1/H_1 \cong G_2/H_2$ such
that  $H$ coincides with the preimage in $G_1\times G_2$ of the
graph of $\gamma$ in $G_1/H_1\times G_2/H_2$.
\end{lem}

\begin{rem}
\label{rg}
\begin{enumerate}
\item[(i)]
If $H_1=G_1, H_2=G_2$ then $H=G_1\times G_2$. If $H_1=\{1\},
H_2=\{1\}$ then $G_1\cong G_2 \cong G$.
\item[(ii)]
If $G_1$ and $G_2$ are disjoint finite groups then one may easily
check that every subgroup of  $G_1\times G_2$ that maps
surjectively on each of the factors coincides with  $G_1\times
G_2$.
\item[(iii)]
 If $G_1=G_2=G$ is a finite simple group then one may easily check that
  either
$H_1=G_1, H_2=G_2,H=G_1\times G_2$ or $H_1=\{1\}, H_2=\{1\}$ and
$G_1\cong G_2 \cong G$.
\end{enumerate}
\end{rem}

\begin{prop}
\label{disjoint} Let $K$ be a field of characteristic different
from $2$ and  $K_a$ its algebraic closure. Let $f(x), h(x)\in
K[x]$ -are irreducible polynomials without multiple roots of
degree
 $n\ge 3$ and
$m\ge 3$ respectively. Suppose that the Galois groups  $\Gal(f)$
of $f$ and $\Gal(h)$  of $h$ are disjoint. Suppose that $f(x)$ is
very nice and $h(x)$ is nice.

Then either
$$\Hom(J(C_f), J(C_h))=0, \Hom(J(C_h),J(C_f))=0$$
or $\fchar(K)>0$ and both jacobians $J(C_f)$ and  $J(C_h)$ are
supersingular abelian varieties.
\end{prop}

\begin{proof}[Proof of Proposition \ref{disjoint}]
Let $K(\RR_f)$ and $K(\RR_h)$ be the splitting fields of
 $f$ and $h$ respectively and  let $L$ be the compositum of
 $K(\RR_f)$ and $K(\RR_h)$. Then the Galois group $\Gal(L/K)$
of $L/K$ may be viewed as a certain subgroup of  $\Gal(f)\times
\Gal(h)$ that maps surjectively  (under the projection maps)  on
each of the factors  $\Gal(f)$ and $\Gal(h)$. It follows from
Remark  \ref{rg}(ii)  and disjointness of $\Gal(f)$ and $\Gal(h)$
that
 $\Gal(L/K)$ coincides with the product
$\Gal(f)\times \Gal(h)$. This means that the extensions
  $K(\RR_f)$ and $K(\RR_h)$ are linearly disjoint over
$K$ and Proposition   \ref{disjoint} follows readily from Theorem
\ref{main}.
\end{proof}

\begin{proof}[Proof of Corollary  \ref{NneM}] It follows from Theorem
\ref{endo} that $\End(J(C_f))=\Z$. Therefore if $\Hom(J(C_f),
J(C_h))\ne 0$, then $\dim(J(C_h)) \ge \dim(J(C_f))$. It follows
that  $\deg(h) \ge 5$ but if $\fchar(K)>0$ then $m=\deg(h) \ge 9$.
Applying again Theorem  \ref{endo}, we observe that
$\End(J(C_h))=\Z$. It follows that if $\Hom(J(C_f), J(C_h))\ne 0$
then  $\dim(J(C_h))= \dim(J(C_f))$. The last equality means that
either $n=m$ or $n$ is even and $m=n-1$ or $m$ is even and
$n=m-1$.

Further, replacing in the case
  $n \ne m,\Gal(f)=\Sn$
the field $K$ by the corresponding quadratic or biquadratic
extension, we may assume that either
$$n \ne m, \quad \Gal(f)=\An, \quad \Gal(h)=\A_m$$
or
$$n=m, \quad \Gal(f)=\Sn, \quad \Gal(h)=\A_m=\An.$$
Notice that in both cases the groups   $G_1:=\Gal(f)$ and
$G_2:=\Gal(h)$ are disjoint. One has only to apply Proposition
 \ref{disjoint}.
\end{proof}

In order to prove Corollary \ref{NeqM} we need a certain
elementary assertion from Galois theory. But first let us
introduce the following notations.  Let $L/K$ be the splitting
field of a separable polynomial $f(x) \in K[x]$ of degree $n$.
Then the set of roots $\RR_f$ of $f$ lies in $L$ and generates it
over $K$. This gives rise to the natural embedding
$\Gal(L/K)\hookrightarrow\Perm(\RR_f)$, which we denote by  $r_f$.
On the other hand, every ordering  $\{\alpha_1, \cdots \alpha_n\}$
of elements of  $\RR_f$ (i.e., of roots of $f$) allows us to
identify $\Perm(\RR_f)$ and $\Sn$ and we may view  $r_f$ as
homomorphism
$$r_f:\Gal(L/K)\hookrightarrow\Perm(\RR_f)=\Sn.$$
Notice that for each positive integer $j \le n$ the stabilizer
$\Gal(L/K)_{\alpha_j}$ of $\alpha_j$ in $\Gal(L/K)$ coincides with
the preimage $r_f^{-1}(\Sn^{\{j\}})$ of the subgroup $\Sn^{\{j\}}$
of all permutations that send $j$ into itself.
\begin{lem}
\label{perm} Let us assume that a finite Galois extension
 $L/K$, a positive integer
 $n$ and a transitive permutation group $\Gamma \subset \Sn$ enjoy
 the following properties:

\begin{enumerate}
\item[(i)]
If $\Gal(L/K)$ is the Galois group of $L/K$ then there exists an
embedding $\Gal(L/K)\hookrightarrow \Sn$,  whose image coincides
with $\Gamma$;
\item[(i)]
For each automorphism  $u:\Gamma \to \Gamma$ of $\Gamma$ there is
a permutation  $s\in\Sn$ such that $u(z)=sz s^{-1}\ \forall z \in
\Gamma$.
\end{enumerate}
Suppose that $f(x), h(x) \in K[x]$ are two separable (i.e.,
without multiple roots) irreducible polynomials of degree
 $n$ such that $L$ is a splitting field of each of them.
 Let as assume additionally that there exist orderings
  $\{\alpha_1, \cdots \alpha_n\}$ of roots of
 $f$ and $\{\beta_1, \cdots \beta_n\}$ of roots of
$h$ such that the image of both natural homomorphisms
$$r_f:\Gal(L/K) \hookrightarrow \Perm(\RR_f)=\Sn, \quad r_h:\Gal(L/K) \hookrightarrow
\Perm(\RR_h)=\Sn$$ coincides with  $\Gamma$.

Then if $\alpha$ is a root of $f$, then there exists a root
 $\beta \in L$ of $h$ such that
 $K(\alpha)=K(\beta)$.
\end{lem}

\begin{proof}[Proof of Lemma \ref{perm}]
Clearly,  $\Gal(L/K)\cong \Gamma$ and there is a permutation
$s\in\Sn$ such that
$$r_h(\sigma)=s r_f(\sigma) s^{-1}\quad \forall \sigma \in
\Gal(L/K).$$ If $j=s(i)$ then one may easily check that
$$r_h^{-1}(\Sn^{\{j\}})=r_f^{-1}(\Sn^{\{i\}}),$$ and therefore in
 $\Gal(L/K)$ the stabilizer
$\Gal(L/K)_{\beta_j}$ of $\beta_j$ coincides with the stabilizer
$\Gal(L/K)_{\alpha_j}$ of $\alpha_j$. This means that
$K(\alpha_i)=K(\beta_j)$.
\end{proof}

\begin{proof}[Proof of ]
In light of Corollary \ref{NneM}, we may assume that either
$$\Gal(f)=\Sn, \quad \Gal(h)=\Sn$$
or
$$\Gal(f)=\An, \quad \Gal(h)=\An.$$
Let us assume that the normal field extensions   $K(\RR_f)$ and
$K(\RR_h)$ do {\sl not} coincide (are not isomorphic). Then their
compositum  $L$ coincides neither with $K(\RR_f)$ nor with
 $K(\RR_h)$ and therefore $\Gal(L/K)$ is isomorphic neither to  $\Gal(f)$ nor to $\Gal(h)$.
Applying Remark \ref{rg}(iii) to  $H=\Gal(L/K), G_1=\Gal(f)$ and
$G_2=\Gal(h)$, we conclude that if $\Gal(f)=\An, \Gal(h)=\An$ then
 $\Gal(L/K)=\Gal(f)\times \Gal(h)$, since $H=\Gal(L/K)$ is not
 isomorphic to
 $G_1=\Gal(f)$. Therefore $K(\RR_f)$ and $K(\RR_h)$ are linearly
 disjoint over
$K$ and Corollary \ref{NeqM} follows readily from Theorem
\ref{main}. If $\Gal(f)=\Sn, \Gal(h)=\Sn$ then an easy check up of
the short list of quotients of  $\Sn$ allows us, applying Lemma
 \ref{goursat} to $H=\Gal(L/K),
G_1=\Gal(f)$ and $G_2=\Gal(h)$ , to conclude that either
$\Gal(L/K)=\Gal(f)\times \Gal(h)$ and Corollary \ref{NeqM} follows
readily from Theorem \ref{main} or  $\Gal(L/K)$ contains
$\An\times\An$ and coincides with the following subgroup of index
$2$ in $\Gal(f)\times \Gal(h)=\Sn\times\Sn$:
$$\{(\sigma,\tau) \in \Sn\times\Sn\mid
\sign(\sigma)=\sign(\tau)\}.$$ (Here  $\sign(\sigma)$ is the sign
of $\sigma$.)  Replacing (in the latter case)   $K$ by the
corresponding quadratic extension, we may assume that
$$\Gal(L/K)=\An\times\An,\Gal(f)=\An, \Gal(h)=\An,$$ and the same arguments as in in the previous case prove
Corollary \ref{NeqM}. Therefore in order to prove Corollary
\ref{NeqM}, it suffices to check that the extensions  $K(\RR_f)$
and $K(\RR_h)$ do not coincide. That is what we are going to do
right now.

Let us assume that $K(\RR_f)=K(\RR_h)$. Replacing $K$ by
corresponding quadratic or biquadratic extension, we may assume
that
$$\Gal(f)=\An, \quad \Gal(h)=\An.$$
Let us put $L=K(\RR_f)=K(\RR_h)$. Clearly, $\Gal(L/K)\cong\An$.
Recall that if $n\ge 5$ and $n\ne 6$ then $\Aut(\An)=\Sn$
~\cite[\S 2.17, pp. 299-300]{Suzuki}. Applying Lemma \ref{perm},
we conclude that $K(\alpha)=K(\beta)$ for some roots $\alpha$ of
$f$ and $\beta$ of $h$. However, $K(\alpha)\cong K[x]/fK[x]=K_f$
and $K(\beta)\cong K[x]/hK[x]=K_h$. Therefore the field extensions
$K_f/K$ and $K_h/K$ are isomorphic. Contradiction.
\end{proof}

\section{Examples}
We write $\bar{\Q}$ for the (algebraically closed) field of all
algebraic numbers in   $\C$.

Let us put $f_n(x)=x^n-x-1 \in\Q[x]$ and consider the number field
 $E_n=\Q[x]/f_n\Q[x]$.
 According to Serre ~\cite[remark 2 on p. 45]{SerreG},
 for each positive integer  $n$ the Galois group of
$f_n(x)$ over $\Q$ coincides with $\Sn$. It is also known (ibid)
that for each prime   $p$ the polynomial $\tilde{f}_n(x):=x^n-x-1
\in \F_p[x]$ either has no multiple roots or   $p$ does not divide
$n(n-1)$ and
$$\tilde{f}_n(x)=(x-\frac{n}{1-n})^2 \tilde{w}(x)$$
where $w(x)\in \F_p[x]$ is a polynomial without multiple roots and
 $0\ne \tilde{w}(\frac{n}{1-n})\in \F_p$. Clearly, if
$\tilde{f}_n$ has no multiple roots then, by Hensel's lemma,
$f_n(x)$ splits into a product of linear factors over an
unramified extension of   $\Q_p$ and therefore the field extension
$E_n/\Q$ is unramified over $p$. But if $\tilde{f}_n$ has a
multiple root then the polynomials $\tilde{w}(x)$ and
$(x-\frac{n}{1-n})^2$ are relatively prime in $\F_p(x)$ and,
thanks to well-known generalization of Hensel's lemma ~\cite[\S
3.5, p.105]{J}, $f_n(x)$ splits over $\Q_p$ into a product of a
quadratic polynomial (that is a lifting of  $(x-\frac{n}{1-n})^2$)
and a certain polynomial $w(x)$(that is a lifting of
$\tilde{w}(x)$). In addition, $w(x)$
 splits into a product of linear factors over an
unramified extension of  $\Q_p$. It follows that if  $E_n/\Q$ is
ramified over  $p$ then it does occur exactly at one prime ideal
of the ring of integers of   $E_n$ and the ramification index is
 $2$.

Let us consider the hyperelliptic curve $A_n: y^2=f_n(x)$ defined
over $\Q$ and its jacobian $J(A_n)$. If $n\le 4$ then $J(A_n)$ is
an elliptic curve. If $n \ge 5$ then $\End(J(A_n))=\Z$
\cite{ZarhinMRL}. Therefore the abelian variety $J(A_n)$ is always
absolutely simple. It follows from Corollary \ref{NneM} that if
$n\ge 5, m\ge 3$ and $n \ne m$ then every homomorphism between the
jacobians $J(A_n)$ and $J(A_m)$ defined over $\bar{\Q}$ is zero.
It follows readily that every homomorphism between  $J(A_n)$ and
$J(A_m)$ defined over the field of complex numbers $\C$ is zero.
(Of course, here the only interesting case is when $n=2g+1$ is odd
and $m=2g+2$ is even and absolutely simple abelian varieties
$J(A_{2g+1})$ and $J(A_{2g+2})$ have the same dimension $g$).

According to Schur \cite{Schur}, the Galois group $\Gal(\exp_n)$
of the polynomial
$$\exp_n(x):=1+x+\frac{x^2}{2}+\frac{x^3}{6}+ \cdots +
\frac{x^n}{n!}$$ over $\Q$ is $\Sn$ if $n$ is not divisible by
$4$; if $4\mid n$ then $\Gal(\exp_n)=\An$. Let us consider the
hyperelliptic curve  $B_n: y^2=f_n(x)$ defined over $\Q$ and its
jacobian $J(B_n)$. If $n\le 4$ then $J(B_n)$ is an elliptic curve;
if $n \ge 5$ then $\End(J(B_n))=\Z$ \cite{ZarhinMRL}. Therefore
the abelian variety $J(B_n)$ is always absolutely simple. It
follows from Corollary   \ref{NneM} that if $n\ge 5, m\ge 3$ and
$n \ne m$ then every homomorphism between  $J(B_n)$ and $J(B_m)$
and also between  $J(B_n)$ and $J(A_m)$ defined over $\bar{\Q}$ or
(which is the same)  over $\C$ is zero. It also follows from
Corollary \ref{NneM} that if $n>5$ and $4\mid n$ then every
homomorphism between  $J(B_n)$ and $J(A_n)$ is zero.

Let us prove, using Corollary \ref{NeqM}, that for all  $n>6$
every homomorphism between  $J(B_n)$ and $J(A_n)$ is zero. In
order to do that, let us consider the number field
$H_n=\Q[x]/\exp_n\Q[x]$. Our goal will be reached if we prove that
the fields $E_n$ are $H_n$ non-isomorphic. To that end, using
Chebyshev's theorem (Bertrand's postulate), pick a prime  $p$ with
$$g+1 \le p \le 2g+1$$
where either $n=2g+1$ is odd or  $n=2g+2$ is even. In particular,
$$p \ge g+1 \ge \frac{n}{2}>3.$$

We write
$$\ord_p:\Q^*\to\Z$$
for the discrete valuation of $\Q$ attached to $p$ and normalized
by the condition  $\ord_p(p)=1$ \cite{Koblitz}. One may easily
check that for all positive integers  $i<p$
$$\ord_p\left(\frac{1}{i!}\right)=0,$$
and for all integers $i$ with $p \le i \le n$
$$\ord_p\left(\frac{1}{i!}\right)=-1,$$
except the case
$$n=2g+2=i,p=g+1,\ord_p\left(\frac{1}{i!}\right)=\ord_{g+1}\left(\frac{1}{(2g+2)!}\right)=-2.$$

It follows that the rational number
  $\frac{-1}{p}$ is a {\sl slope} of the
 $p$-adic Newton polygon of
$\exp_n(x)$. The well-known connection between (reciprocal) roots
of a polynomial and slopes of its Newton polygon conclude
 \cite{Koblitz} allows us to conclude that there is a prime ideal
 in the ring of integers of
$H_n$ that divides  $p$ and whose ramification index in $H_n/\Q$
is divisible by $p>3$. Since all the ramification indices in
 $E_n/\Q$ do not exceed
$2$ (see the beginning of this Section), the fields $E_n$ and
$H_n$ are non-isomorphic. Applying Corollary \ref{NeqM} to
$f=f_n,h=\exp_n$, we conclude that for all  $n>6$ every
homomorphism between $J(B_n)$ and $J(A_n)$ defined over $\bar{\Q}$
or (which is the same)  over $\C$ is zero.

Let us turn now to a completely different class of examples. Let
 $p$ be an odd prime,  $k_p$ an algebraically closed field of
 characteristic
 $p$, K=K(t) the field of rational functions in independent variable $t$ with coefficients in
 $k_p$ and
$K_a$ an algebraic closure of $K$. Let an integer $q>1$ be an
integral power of  $p$. Let $d>1$ be a positive integer. Let us
put
$$n=\frac{q^d-1}{q-1}=\#(\P^{d-1}(\F_q))$$
and consider the polynomials
$$f(x)=x^n+tx+1 \in K[x], \quad h(x)=x^n+x+t \in K[x].$$
 According to Abhyankar \cite{Ab}, there are bijections
 $$\RR_{f} \cong \P^{d-1}(\F_q), \quad \RR_{f} \cong
 \P^{d-1}(\F_q),$$
 such that $\Gal(f)$ becomes $\PSL(d,q)$ and
 $\Gal(h)$ becomes $\PGL(d,q)$. Let us assume, in addition, that
 $m>2$. Then both $f(x)$ and $h(x)$ are very nice.
Suppose also that  $d$ and $q-1$ are {\sl not} relatively prime.
Then $\Gal(f)=\PSL(d,q)$ and $\Gal(f)=\PGL(d,q)$ are disjoint (see
 \ref{exd}(iv)).  According to Proposition \ref{disjoint},
 if $J(C_f)$ and $J(C_h)$ are the jacobians of the hyperelliptic
 curves
$$C_f:y^2=f(x), \quad C_h:y^2=h(x)$$
then either both jacobians are supersingular or every homomorphism
between $J(C_f)$ and $J(C_h)$ defined over $K_a$ is zero. However,
by theorem  2.4(iv) of \cite{ZarhinMMJ}, if $(q,d)\ne (3,4)$ then
$$\End(J(C_f))=\Z, \quad \End(J(C_h))=\Z,$$
and therefore both jacobians are not supersingular. Therefore if
$(q,d)\ne (3,4)$ then every homomorphism between  $J(C_f)$ and
$J(C_h)$ defined over $K_a$ is zero.

\vskip  .1cm

 \noindent {\small
Pennsylvania State University, Department of Mathematics,
University Park, PA 16802, USA}

\vskip .1cm

\noindent {\small Institute of Mathematical Problems in Biology,
Russian Academy of Sciences, Pushchino, Moscow Region, Russia  }

 \vskip  .1cm

\noindent {\small {\em E-mail address}: zarhin@math.psu.edu}

\end{document}